\documentclass[11pt,reqno,a4paper]{amsart}
\usepackage{amsmath}
\usepackage{amssymb}
\usepackage{graphicx}
\usepackage{enumerate}
\newcommand{\eps}{\varepsilon}
\newcommand{\N}{\mathbb{N}}

\newcommand{\R}{\mathbb{R}}
\newcommand{\T}{\mathbb{T}}

\renewcommand{\P}{\mathbb{P}}

\DeclareMathOperator{\ve}{\varepsilon}
\DeclareMathOperator{\p}{\mathbb{P}}
\DeclareMathOperator{\om}{\omega}

\title{Recurrence for random dynamical systems}
\date{June 25, 2009}
\author{Philippe Marie and Jerome Rousseau}
\address{Philippe Marie, Universit\'e d'Aix-Marseille, Centre de physique th\'eorique, UMR 6206 CNRS, Campus de Luminy, Case 907 - 13288 Marseille cedex 9 and Universit\'e du Sud, Toulon-Var, France.}
\email{pmarie@cpt.univ-mrs.fr}
\urladdr{http://www.cpt.univ-mrs.fr/~marie}
\address{Jerome Rousseau, Universit\'e Europ\'eenne de Bretagne, Universit\'e de Brest, Laboratoire de Math\'ematiques CNRS UMR 6205, 6 avenue Victor le Gorgeu, CS93837, F-29238 Brest Cedex 3, France}
\email{jerome.rousseau@univ-brest.fr}
\urladdr{http://pageperso.univ-brest.fr/~rousseau}
\begin{document}
\newtheorem{lem}{Lemma}
\newtheorem{theorem}[lem]{Theorem}
\newtheorem{defi}[lem]{Definition}
\newtheorem{prop}[lem]{Proposition}
\newtheorem{corollary}[lem]{Corollary}
\newtheorem{example}[lem]{Example}
\newtheorem*{remark}{Remark}
\newtheorem{definition}[lem]{Definition}

\maketitle
\begin{abstract}
This paper is a first step in the study of the recurrence behavior in random dynamical systems and randomly perturbed dynamical systems. In particular we define a concept of quenched and annealed return times for systems generated by the composition of random maps. We moreover prove that for super-polynomially mixing systems, the random recurrence rate is equal to the local dimension of the stationary measure.
\end{abstract}

\section{Introduction}
 Random dynamical systems (e.g. \cite{K,KL,L,O}) and  quantitative description of the recurrence in deterministic dynamics  (e.g. \cite{B,BS,S}) have been deeply studied in the last years. But despite of the success of these both surveys, we did not find in the literature any result about the recurrence of random dynamical systems. The purpose of this paper is to define the main needed objects to start this study and to give  first results to describe the recurrence behavior in the random case.

The evolution of a random dynamical system generated by \textit{i.i.d.} random transformations will be as described in what follows: we consider an indexed family of transformations $\{T_{\lambda}\}_{\lambda\in\Lambda}$ defined on a compact Riemannian manifold $X$ and a probability measure $\mathcal{P}$ on a metric space $\Lambda$.  Let $(\lambda_n)_{n\geq 0}$ be an \textit{i.i.d.} stochastic process with common distribution $\mathcal{P}$, a random evolution of an initial state $x\in X$  will be generated by a realisation of this process, say $\underline{\lambda}:=(\lambda_0,\lambda_1,\ldots,\lambda_n,\ldots)\in\Lambda^{\mathbb{N}}$, such as for every $n\geq 0$:
$$
x_n^{\underline{\lambda}}:=T_{\lambda_{n-1}}\circ\ldots T_{\lambda_0}x
$$ 
The set $\{x_n^{\underline{\lambda}}\}_{n=0}^{\infty}\subset X$ will be the random orbit of $x$ associated to the realisation $\underline{\lambda}\in\Lambda^{\mathbb{N}}$. 
We actually do not require for the process $(\lambda_n)_{n\geq0}$ to be independent. We consider a weaker assumption that we explain below which allow us to get results available for more general random dynamical systems. In fact, we focus on systems for which there exists  a stationary measure; in the \textit{i.i.d.} case a stationary measure is a probability measure $\mu$ on $X$ such that :
$$
\mu(A)=\int_{\Lambda}\mu(T_{\lambda}^{-1}A)d\mathcal{P}(\lambda)
$$
for any Borelian $A$ of $X$.
The probability measure $\mu$ over the phase space $X$ will play the role of the invariant one in the random setting. It is natural to wonder if the classical results of deterministic recurrence can be applied to the random orbits evolution. There are two ways of studying the long-term behavior of random orbits; we can  describe them with respect to the stationary measure, firstly, fixing a realisation of the random component  and secondly, taking into account the noise globally. The first case is a statistical study only over the space $X$, while the second one is a statistical study over both the space $X$ and the noise space $\Omega$. The first point of view is usually called \textit{quenched} and the second one \textit{annealed}. The annealed case provides  a better global understanding of random dynamical systems as ergodic theory provides a great scope of  deterministic orbits description.

First we will introduce some new definitions that are the counterpart of usual recurrence objects in the random framework. More precisely, for random dynamical systems generated by the composition of random transformations we define a \textit{quenched} and an \textit{annealed} version of return times. The \textit{annealed }one is defined in the i.i.d. case by
\[\mathbf{T}_B(x):=\int_{\Lambda^\N}\inf\{n>0:T_{\lambda_{n-1}}\circ\ldots\circ T_{\lambda_0}x\in B\}d\mathcal{P}^{\mathbb{N}}(\underline\lambda).\]

As we said above, this quantity appears to us as the best quantity for a relevant   quantitative analysis of the recurrence in random dynamical systems. We will justify the choice of this definition applying classical results of deterministic recurrence   (Poincar\'e's theorem, Kac's lemma) to a usual representation of random maps, namely the skew product transformation, and a good Borelian subset. We will moreover see that the \textit{annealed} return time appears naturally in this study.

Then, in section 3, using a recent result from Rousseau and Saussol \cite{RS} about the recurrence for observations applied to the skew product representation, we establish a link between random return times and local dimension of stationary measures in both cases: \textit{quenched} and \textit{annealed}. Moreover, in our main theorem (Theorem~\ref{thegal}),we prove that, for systems whose decay of correlations hold super-polynomially fast for Lipschitz observables, the random recurrence rates are equal to the local dimension of the stationary measure. This is a first step in the description of the recurrence for non-deterministic dynamical systems.

In the fourth section we will give two examples of random toral automorphisms for which our results hold. The first one is a non-i.i.d.\ random dynamical system where the random maps are themselves chosen according to a Markov process. The second one deals with i.i.d.\ hyperbolic toral automorphisms.

We will then  consider the particular case where random dynamical systems are used as a small stochastic perturbation of an initial dynamics. We will apply our results to two examples studied by Baladi and Young \cite{BY}: expanding maps of the circle and piecewise expanding maps of the interval.

The last part is devoted to the proofs of our main theorems.
\section{Random dynamical systems and return time}
Let us first recall some generalities about random dynamical systems. Let $\Omega$ be a metric space with $\mathcal{B}(\Omega)$ its Borelian $\sigma$-algebra and $(\Omega,\mathcal{B}(\Omega),\mathbb{P})$ a probability space which represent the space of the randomness. Let $\vartheta$ be a $\P$-preserving map which represents the time evolution of the randomness. $\Omega$ indexes a family of maps $\{T_{\om}\}_{\om\in\Omega}$ from a compact Riemannian manifold $X$ into itself. A random dynamical system $\mathcal{T}$ on $X$ over $(\Omega,\mathcal{B}(\Omega),\mathbb{P},\vartheta)$ is generated by mappings $T_\omega$ so that the map $(\omega,x)\rightarrow T_\omega x$ is measurable and so that the following cocycle property holds:
$$
T^{n+m}_{\omega}=T^n_{\vartheta^{m}\omega}\circ T^m_{\omega}
$$
Each  realisation  $\omega\in\Omega$ of the random component of the system defines a random evolution law, under the action of which, the evolution of an initial state $x\in X$ will be given by its random orbit, defined by the set $\{x_n^{\omega}\}$, where for any $n\in\mathbb{N}$: 
$$
x_n^{\omega}:=T_{\vartheta^{n-1}\om}\circ\ldots\circ T_{\om}x:=T^n_\omega x.
$$
One useful representation of this system is given by the following skew product transformation:
\begin{align*}
S:~~\Omega\times X &\longrightarrow \Omega\times X\\
(\omega,x) &\longmapsto (\vartheta\omega,T_{\om}x).
\end{align*}
It is straightforward to see that for a fixed $(\omega,x)\in\Omega\times X$ the iteration of this map will generate on the second component the random orbit of $x$ associated to the realisation $\omega$ since
$S^n(\omega,x)=(\vartheta^n\omega,T_{\omega}^nx)$.
This drives to consider invariant measures of the skew product map to study the behavior of the random dynamical system. Let us suppose that $\nu$ is an invariant probability measure for $S$ and define the canonical projection over $\Omega$:
\[
\begin{aligned}
\pi_{\Omega}:\Omega\times X &\rightarrow \Omega\\
(\omega,x) &\mapsto \omega
\end{aligned}
\]
Then, denoting by $F_*\mu$ the pushforward measure of $\mu$ by a measurable map $F$ \footnote{\textit{i.e.} for all measurable subset $A$, $F_*\mu(A)=\mu(F^{-1}A)$},  $S_*\nu=\nu$ implies $\vartheta_*((\pi_{\Omega})_*\nu)=(\pi_{\Omega})_*\nu$ and thus the marginal measure of $\nu$ over $(\Omega,\mathcal{B}(\Omega))$ must be invariant for the dynamics $\vartheta$. Since the random component of the system is considered to be given a priori \cite{Arnold}, we will focus on invariant measures     satisfying $(\pi_{\Omega})_*\nu=\P$.
\begin{definition}
A probability measure $\nu$ is invariant for the random dynamical system $\mathcal{T}$ if:
\begin{enumerate}[(i)]
\item $S_*\nu=\nu$
\item $(\pi_{\Omega})_*\nu=\p$.
\end{enumerate}
\end{definition}

From now on,  we  deal with random dynamical systems for which the skew product invariant measure is a product measure $\nu=\P\otimes\mu$, this case  is in fact quite general since it includes the important case we considered in introduction where the maps  are chosen independently with the same distribution $\mathcal{P}$. We recover this situation setting  $\P=\mathcal{P}^\N$,  $\vartheta=\sigma$, the left shift defined on $\Omega=\Lambda^{\mathbb{N}}$ and $T_\omega=T_{\lambda_0}$ for $\omega=(\lambda_0,\lambda_1,...,\lambda_n,...)$ in the skew product definition.
\begin{definition}
We say that $\mu$ is a stationary measure for the random dynamical system if the measure $\nu=\mu\otimes\P$ is invariant for the skew product.
\end{definition}
This vocabulary comes from the fact that i.i.d.\ random maps generate Markov processes, namely, a random orbit is a Markov process whose transition probabilities are given by:
\begin{equation} \label{Markov}
P(x,B)=\P\{\omega\in\Omega~:~T_{\omega}x\in B\}  
\end{equation}
A huge literature is devoted to the study of this special kind of Markov processes, see for example \cite{BL, BM, DF}. This property of random orbits provides a great interest to this study from the probabilistic point of view. 
 We remark that for i.i.d.\ random dynamical systems under weak assumptions such measures always exist (for example if $T_\omega$ is continuous for all $\omega\in\Omega$ \cite{K}). 
\begin{prop}
If $\mathcal{T}=(M,\{T_{\lambda}\}_{\lambda\in\Lambda}, \mathcal{P},\sigma)$ is an i.i.d.\ random dynamical system, then a probability measure $\mu$ is stationary if and only if for all Borelian $B$ it verifies (see \cite{Arnold} for instance):
$$
\int_{\Lambda}(T_{\lambda})_*\mu(B)d\mathcal{P}(\lambda)=\mu(B)
$$
\end{prop} 
That is to say that a stationary measure is a measure which is invariant in average over the random component.

Therefore, the study of the deterministic dynamical system $(\Omega\times X, S, \P\otimes\mu)$ will provide lots of informations about the behavior of random orbits.

We now introduce a return time concept for random orbits.
Since $(\Omega\times X,S,\p\otimes \mu)$ is a deterministic dynamical system, we can define the first return time into a Borelian subset $A\times B\subset \Omega\times X$ as usual for $(\om,x)\in A\times B$:
\[\tau_{A\times B}^S(\omega,x)=\inf\{k> 0~:~S^k(\omega,x)\in A\times B\}.\]
Choosing $A=\Omega$ we remark that $S^n(\omega,x)\in \Omega\times B$ if and only if $T_{\omega}^nx\in B$. This naturally drives to the following definition:
\begin{definition} For a fixed $\omega\in \Omega$ the first quenched random return time in a measurable subset $B\subset X$ of the random orbit starting from a point $x\in B$ is:
\begin{eqnarray*}
\tau_B^{\omega}(x)&=&\tau_{\Omega\times B}^S(\omega,x)\\
&=&\inf\{k>0\,:\,T_{\vartheta^{k-1}\omega}\circ\dots\circ T_{\omega}x\in B\}.
\end{eqnarray*}
\end{definition}
Remark that  the Poincar\'e recurrence theorem applied to the skew product ensures that this quantity is almost everywhere finite.

From now on we assume that $X$ is a metric space with metric $d$. For $\omega\in\Omega$, we are interested in the behavior as $r\rightarrow0$ of the quenched random return time of a point $x\in X$ into the open ball $B(x,r)$ defined by
\begin{eqnarray*}\tau_{r}^{\omega}(x)&:=&\inf\left\{k>0:T_{\vartheta^{k-1}\omega}\circ...\circ T_{\omega}x\in B(x,r)\right\}\\&=&\inf \left\{k>0: d(T_{\vartheta^{k-1}\omega}\circ...\circ T_{\omega}x,x)<r\right\}.\end{eqnarray*}
\begin{definition}
The random dynamical system $\mathcal{T}$ on $X$ over $(\Omega,\mathcal{B}(\Omega),\mathbb{P},\vartheta)$ with a stationary measure $\mu$ is called random-aperiodic if 
\[\P\otimes\mu\{(\omega,x)\in\Omega\times X:\exists n\in \mathbb{N},\,T_{\vartheta^{n-1}\omega}\circ\ldots\circ T_{\omega}x=x\}=0.\]
\end{definition}
Let $T_0=\mathcal{R}_\alpha$ an irrational rotation of the circle for an irrational number $\alpha$ and $T_1$ the identity map of the circle. The i.i.d.\ random dynamical system constructed with this two maps chosen with the same probability $P(0)=P(1)=\frac12$ is not random-aperiodic and we have $\P\{\om\in\Omega~: ~\tau^\omega_r(x)=1\}=P(1)=\frac{1}{2}$ for all $r>0$. This means that after one iteration half of the points did not rotate and $\left(\frac{1}{2}\right)^n$ of them after $n$ iterations. Anyway,  almost every point will eventually rotate and their dynamics will be quite interesting providing that we wait for an enough  long time. To avoid this kind of problem with the first return time for non-random-aperiodic system we need to introduce the non-instantaneous return times (more details are presented on the Section 2.3 of \cite{RS}). 
\begin{definition}
Let $r>0$. For $x\in X$, $\omega\in\Omega$ and $p\in\mathbb{N}$ we define the $p$-non-instantaneous quenched random return time:
\[\tau_{r,p}^{\omega}(x):=\inf \left\{k>p~:~ d(T_{\vartheta^{k-1}\omega}\circ...\circ T_{\omega}x,x)<r\right\}.\]
Then we define the non-instantaneous quenched random lower and upper recurrence rates:
\[\underline{R}^{\omega}(x):=\lim_{p\rightarrow\infty}\liminf_{r\rightarrow0}\frac{\log\tau_{r,p}^{\omega}(x)}{-\log r}\qquad\textrm{and}\qquad\overline{R}^{\omega}(x):=\lim_{p\rightarrow\infty}\limsup_{r\rightarrow0}\frac{\log\tau_{r,p}^{\omega}(x)}{-\log r}.\]
\end{definition}

Finally we recall that \emph{the lower and upper pointwise or local dimension} of a Borel probability measure $\mu$ on X at a point $x\in X$ are defined by
\[\underline{d}_\mu(x)=\underset{r\rightarrow0}{\liminf}\frac{\log\mu\left(B\left(x,r\right)\right)}{\log r}\qquad\textrm{and}\qquad\overline{d}_\mu(x)=\underset{r\rightarrow0}{\limsup}\frac{\log\mu\left(B\left(x,r\right)\right)}{\log r}.\]
These dimensions are linked with the Hausdorff dimension of the measure (see e.g. \cite{F} for more details on dimensions):
\begin{remark}If $X\subset\R^n$ for some $n\in\N^*$, for a finite Borel measure $\mu$ on $X$
\[\dim_H\mu=\textrm{ess-sup }\underline{d}_\mu.\]
\end{remark}


\section{Quenched recurrence rate, annealed recurrence rate and pointwise dimension}


We now state our results about the recurrence behavior for the random maps. We first focus on the quenched case (random orbits) in section 3.1 and then on the annealed one (random recurrence) in section 3.2.

\subsection{Recurrence for random orbits}
\begin{theorem}\label{th1}Let $\mathcal{T}$ be a random dynamical system on $X$ over $(\Omega,\mathcal{B}(\Omega),\mathbb{P},\vartheta)$ with a stationary measure $\mu$. For $\P\otimes\mu$-almost every $(\omega,x)\in\Omega\times X$
\[\underline{R}^{\omega}(x)\leq\underline{d}_\mu(x)\qquad\textrm{and}\qquad\overline{R}^{\omega}(x)\leq\overline{d}_\mu(x).\]
\end{theorem}
\begin{proof}This theorem is proved using Theorem 2 of \cite{RS} applied to $(\Omega\times X,\mathcal{B}(\Omega\times X),\P\otimes\mu ,S)$ with the observation $f$ defined by
\begin{eqnarray*}
f\,\,:\Omega\times X&\longrightarrow& X\\
(\omega,x)&\longmapsto& x.
\end{eqnarray*}
With this observation, for all $(\omega,x)\in\Omega\times X$, for all $r>0$ and for all $p\in\N$ we identify the return time for the observation
\[\tau^f_{r,p}(\omega,x)=\tau^\omega_{r,p}(x),\] 
the pushforward measure
\[f_*(\P\otimes\mu)=\mu,\]
and the pointwise dimensions for the observation
\[\underline{d}^f_{\P\otimes\mu}(x)=\underline{d}_\mu(x)\qquad\textrm{and}\qquad \overline{d}^f_{\P\otimes\mu}(x)=\overline{d}_\mu(x).\]
\end{proof}
This theorem is also satisfied for the first return time:
\begin{corollary}Let $\mathcal{T}$ be a random dynamical system on $X$ over $(\Omega,\mathcal{B}(\Omega),\mathbb{P},\vartheta)$ with a stationary measure $\mu$. For $\P\otimes\mu$-almost every $(\omega,x)\in\Omega\times X$
\[\liminf_{r\rightarrow0}\frac{\log\tau_{r}^{\omega}(x)}{-\log r}\leq\underline{d}_\mu(x)\qquad\textrm{and}\qquad\limsup_{r\rightarrow0}\frac{\log\tau_{r}^{\omega}(x)}{-\log r}(x)\leq\overline{d}_\mu(x).\]
\end{corollary}
\begin{proof} Since for all $(\omega,x)\in\Omega\times X$, for all $r>0$ and for all $p\in\N$ we have $\tau_{r}^{\omega}(x)\leq\tau_{r,p}^{\omega}(x)$ then
\[\liminf_{r\rightarrow0}\frac{\log\tau_{r}^{\omega}(x)}{-\log r}\leq\underline{R}^{\omega}(x)\qquad\textrm{and}\qquad\limsup_{r\rightarrow0}\frac{\log\tau_{r}^{\omega}(x)}{-\log r}(x)\leq\overline{R}^{\omega}(x)\]
and the corollary is proved by Theorem~\ref{th1}.
\end{proof}
Even if the inequalities in Theorem~\ref{th1} can be strict, with more assumptions on the random dynamical system one can prove that the equalities hold. This drives us to introduce the decay of correlations for a random dynamical system:
\begin{definition}
A random dynamical system with a stationary measure $\mu$ has a super-polynomial decay of correlations if for all $n\in\mathbb{N}^{*}$ and $\psi$ , $\varphi$ Lipschitz observables from $X$ to $\mathbb{R}$:
$$
\vert\int_X \int_{\Omega}\psi(T_{\vartheta^{n-1}\om}\circ\ldots\circ T_{\om}x)\varphi(x)d\P(\omega)d\mu(x)-\int_X\psi d\mu \int_X \varphi d\mu\vert\leq \Vert \psi\Vert\Vert \varphi\Vert\theta_n
$$
with $\lim_{n\rightarrow\infty}\theta_nn^p=0$ for any $p>0$ and where $\Vert.\Vert$ is the Lipschitz norm.
\end{definition}

\begin{theorem}\label{spdoc} Let $\mathcal{T}$ be a random dynamical system on $X$ over $(\Omega,\mathcal{B}(\Omega),\mathbb{P},\vartheta)$ with a stationary measure $\mu$. If the random dynamical system has a super-polynomial decay of correlations then
\[\underline{R}^{\omega}(x)=\underline{d}_\mu(x)\qquad\textrm{and}\qquad\overline{R}^{\omega}(x)=\overline{d}_\mu(x)\]
for $\P\otimes\mu$-almost every $(\omega,x)\in\Omega\times X$ such that $\underline{d}_\mu(x)>0$.
\end{theorem}
\begin{proof}As previously, this theorem is proved using Theorem 5 of \cite{RS} applied to $(\Omega\times X,\mathcal{B}(\Omega\times X),\P\otimes\mu ,S)$ with the observation $f$ defined by
\begin{eqnarray*}
f\,\,:\Omega\times X&\longrightarrow& X\\
(\omega,x)&\longmapsto& x.
\end{eqnarray*}
We emphasize that in Theorem 5 of \cite{RS} the condition of the super-polynomial decay of correlations is applied to the skew product, in fact in our case this is not necessary. A weaker assumption, the super-polynomial decay of correlations for the random dynamical system, allows us to prove Lemma 15 of \cite{RS} for the observation $f$ defined above. The proof of this lemma is identical since in our case the decay of correlations is used in a function which does not depend on $\omega$.
\end{proof}
The following proposition establishes that the non-instantaneous return time notions are not necessary for random-aperiodic systems.
\begin{prop}\label{ra} If the system is random-aperiodic, for $\P\otimes\mu$-almost every $(\omega,x)\in\Omega\times X$ and for every $p\in\N$, $\tau^{\omega}_{r,p}(x)=\tau^{\omega}_r(x)$ for $r$ small enough. 
\end{prop}
\begin{proof}
Since the system is random-aperiodic, for $\P\times\mu$-almost all $(\omega,x)\in\Omega\times X$, for all $k>0$, $d(T_{\vartheta^{k-1}\omega}\circ\ldots\circ T_{\omega}x,x)>0$, then $\tau^{\omega}_r(x)\underset{r\rightarrow0}\rightarrow+\infty$. Thus, for all $p>0$, there exists $r(p,x,\omega)$ such that for every $r< r(p,\omega,x)$, we have $\tau_{r}^{\omega}(x)>p$. Therefore, for every $r< r(p,\omega,x)$, we have $\tau_{r}^{\omega}(x)=\tau_{r,p}^{\omega}(x)$.
\end{proof}
The previous theorem is proved for the instantaneous quenched random lower and upper recurrence rate:
\begin{corollary}\label{cospdoc} Let $\mathcal{T}$ be a random dynamical system on $X$ over $(\Omega,\mathcal{B}(\Omega),\mathbb{P},\vartheta)$ with a stationary measure $\mu$. If the random dynamical system is random-aperiodic and has a super-polynomial decay of correlations then
\[\liminf_{r\rightarrow0}\frac{\log\tau_{r}^{\omega}(x)}{-\log r}=\underline{d}_\mu(x)\qquad\textrm{and}\qquad\limsup_{r\rightarrow0}\frac{\log\tau_{r}^{\omega}(x)}{-\log r}=\overline{d}_\mu(x)\]
for $\P\otimes\mu$-almost every $(\omega,x)\in\Omega\times X$ such that $\underline{d}_\mu(x)>0$.
\end{corollary}
\begin{proof}This is just a consequence of Theorem~\ref{spdoc} and Proposition~\ref{ra}.
\end{proof}


\subsection{Random recurrence}
For $B$ a Borelian subset of $X$, Poincar\'e's recurrence theorem applied to $S$ ensures that the first return time is finite for $\P\otimes \mu$-almost all $(\omega,x)\in \Omega\times B$. The other main result in deterministic recurrence is Kac's lemma. To study what happens in the random case, we first need to recall the concept of ergodicity in the random setting:
\begin{definition}
 We say that the random dynamical system $\mathcal{T}$ on $X$ over $(\Omega,\mathcal{B}(\Omega),\mathbb{P},\vartheta)$ is ergodic with respect to a stationary measure $\mu$ if the deterministic system $(\Omega\times X,S,\P\otimes\mu)$ is ergodic.
\end{definition}
Let us now suppose that the random dynamical system is ergodic then we can apply Kac's lemma and we get:
\begin{eqnarray}
\int_{\Omega\times B}\tau_{\Omega\times B}^S(\omega,x)\:d\P\otimes \mu(\omega,x)&=&\int_{\Omega\times B}\tau_B^{\omega}(x)   \:d\P\otimes \mu(\omega,x)\nonumber\\
&=&\int_B \int_{\Omega}\tau_B^{\omega}(x)  \:d\P(\omega)d\mu(x)\nonumber\\
&=1.\label{kac}
\end{eqnarray}
Which is an argument to study the quantity: $$\int_{\Omega}\tau_B^{\omega}(x)  \:d\P(\omega).$$
Therefore, the latter naturally appears in the study of random recurrence and while the quenched version of return times will allow a description realisation by realisation, this one will provide a global description of the system, including the random part in its totality. This is coherent with the ergodic theory point of view and the classical statistical study of dynamical systems. Let us moreover remark that averaged quantities are also omnipresent in the study of i.i.d.\ random dynamical systems, see for example the Viana's course \cite{Vi}, where the random evolution operator, respectively the transfert one, is the average of the evolution, respectively the transfert, operators associated to each random maps $T_{\omega}$ over the random component. We have already remarked in the first section that stationary measures were defined by an average too.
\begin{definition}
The first annealed random return time for $x\in B$ in the set $B$ is:
$$
\mathbf{T}_B(x)=\int_{\Omega}\tau_B^{\omega}(x)  \:d\P(\omega).
$$
\end{definition}
With this definition and \eqref{kac}, we immediately get a random version of Kac's lemma:
\begin{prop}[Random Kac's lemma]
Let $\mathcal{T}$ be an ergodic random dynamical system with respect to the stationary probability measure $\mu$, then for any Borelian subset $B$ of $X$ such that $\mu(B)>0$, we have:
\[\int_B \mathbf{T}_B(x)d\mu(x)=1.\]
\end{prop}

As previously, we are interested in the behaviour of the return time of $x\in X$ into $B(x,r)$ when $r\rightarrow0$, which drives to the following definitions:
\begin{definition}
Let $r>0$ and $p\in\N$. We define the $p$-non-instantaneous annealed random return time
\[\mathbf{T}_{r,p}(x):=\int_{\Omega}\tau_{r,p}^{\omega}(x)d\P(\omega).\]
When $p=0$, we denote this return time $\mathbf{T}_{r}(x)$. Then we define the non-instantaneous annealed random lower and upper recurrence rates
\[\underline{\mathbf{R}}(x):=\lim_{p\rightarrow\infty}\liminf_{r\rightarrow0}\frac{\log\mathbf{T}_{r,p}(x)}{-\log r}\qquad\overline{\mathbf{R}}(x):=\lim_{p\rightarrow\infty}\limsup_{r\rightarrow0}\frac{\log\mathbf{T}_{r,p}(x)}{-\log r}.\]
\end{definition}
The main results are the following two theorems (whose proofs can be found in Section~\ref{proof}):
\begin{theorem}\label{thborne}Let $\mathcal{T}$ be a random dynamical system on $X$ over $(\Omega,\mathcal{B}(\Omega),\mathbb{P},\vartheta)$ with a stationary measure $\mu$. For $\mu$-almost every $x\in X$
\[\underline{\mathbf{R}}(x)\leq\underline{d}_\mu(x)\qquad\textrm{and}\qquad\overline{\mathbf{R}}(x)\leq\overline{d}_\mu(x).\]
\end{theorem}
\begin{theorem}\label{thegal}Let $\mathcal{T}$ be a random dynamical system on $X$ over $(\Omega,\mathcal{B}(\Omega),\mathbb{P},\vartheta)$ with a stationary measure $\mu$. If the random dynamical system has a super-polynomial decay of correlations then
\[\underline{\mathbf{R}}(x)=\underline{d}_\mu(x)\qquad\textrm{and}\qquad\overline{\mathbf{R}}(x)=\overline{d}_\mu(x)\]
for $\mu$-almost every $x$ such that $\underline{d}_\mu(x)>0$. 
\end{theorem}
The first inequality is still satisfied using the first return time:
\begin{corollary}\label{coborne}Let $\mathcal{T}$ be a random dynamical system on $X$ over $(\Omega,\mathcal{B}(\Omega),\mathbb{P},\vartheta)$ with a stationary measure $\mu$. For $\mu$-almost every $x\in X$
\[\liminf_{r\rightarrow0}\frac{\log\mathbf{T}_r(x)}{\log r}\leq\underline{d}_\mu(x)\qquad\textrm{and}\qquad\limsup_{r\rightarrow0}\frac{\log\mathbf{T}_r(x)}{\log r}\leq\overline{d}_\mu(x).\]
\end{corollary}
\begin{proof} Since for every $x\in X$ and for every $p\in\N$, $\mathbf{T}_r(x) \leq\mathbf{T}_{r,p}(x)$, the result follow from Theorem~\ref{thborne}.
\end{proof}
As in the previous section, for random-aperiodic dynamical systems we do not need the non-instantaneous return time (see Section~\ref{proof} for the proof):
\begin{prop}\label{propegal}
Let $\mathcal{T}$ be a random dynamical system on $X$ over $(\Omega,\mathcal{B}(\Omega),\mathbb{P},\vartheta)$ with a stationary measure $\mu$. If the random dynamical system is random-aperiodic and has a super-polynomial decay of correlations then
\[\liminf_{r\rightarrow0}\frac{\log\mathbf{T}_r(x)}{\log r}=\underline{d}_\mu(x)\qquad\textrm{and}\qquad\limsup_{r\rightarrow0}\frac{\log\mathbf{T}_r(x)}{\log r}=\overline{d}_\mu(x)\]
for $\mu$-almost every $x$ such that $\underline{d}_\mu(x)>0$. 
\end{prop}

\section{Random toral automorphisms}
We now give two examples of random dynamical systems for which our theorems hold. We emphasize that the first example is a non-i.i.d.\ random dynamical system where there exists a stationary measure.
\subsection{Non-i.i.d.\ random linear maps}
Let $X=\T^1$ be the one-dimensional torus. Consider the two linear maps which preserve Lebesgue measure $Leb$ on $X$
\[\begin{array}{rcccrcc}
T_1:X&\longrightarrow& X & \textrm{ and }&T_2:X&\longrightarrow& X\\
x&\longmapsto& 2x& & x&\longmapsto& 3x.
\end{array}\]
The random orbit is constructed by choosing one of these two maps following a Markov process with the stochastic matrix
\[A=\begin{pmatrix}1/2&1/2\\1/3&2/3\end{pmatrix}.\]
In fact, this random dynamical system is represented by the following skew product
\begin{eqnarray*}
S:\Omega\times X&\longrightarrow& \Omega\times X\\
(\omega,x)&\longmapsto&(\vartheta(\omega),T_\omega x)
\end{eqnarray*}
with $\Omega=[0,1]$, $T_\omega=T_1$ if $\omega\in[0,2/5)$, $T_\omega=T_2$ if $\omega\in[2/5,1]$ and where $\vartheta$ is the following piecewise linear map
\[\vartheta(\omega)=\left\{\begin{array}{lll}
2\omega\qquad&\textrm{if}&\omega\in[0,1/5)\\
3\omega-1/5\qquad&\textrm{if}&\omega\in[1/5,2/5)\\
2\omega-4/5\qquad&\textrm{if}&\omega\in[2/5,3/5)\\
3\omega/2-1/2 \qquad&\textrm{if}&\omega\in[3/5,1].\\
\end{array}\right.\]
Since $\vartheta$ preserves Lebesgue measure $\P=Leb$, the skew product $S$ is $Leb\otimes Leb$-invariant and so $Leb$ is a stationary measure for this random dynamical system. One can easily see that this system is random-aperiodic. Since $S$ has an exponential decay of correlations \cite{Ba}, Corollary~\ref{cospdoc} and Proposition~\ref{propegal} hold, so for $Leb\otimes Leb$-almost every $(\omega,x)\in[0,1]\times\T^1$
\[\lim_{r\rightarrow0}\frac{\log \tau_{r}^{\omega}(x)}{-\log r}=1\]
and for $Leb$-almost every $x\in\T^1$
\[\lim_{r\rightarrow0}\frac{\log\left[\mathbf{T}_{r}(x)\right]}{-\log r}=1.\]

\subsection{Random hyperbolic toral automorphisms}
Let $X=\mathbb{T}^2$, we recall that a hyperbolic toral automorphism is a map $A:\mathbb{T}^2\rightarrow\mathbb{T}^2$ acting through the matrix $x\mapsto Ax$ (mod 1), such that the matrix $A$ has integer entries, eigenvalues with modulus different from 1 and  $\det A=\pm 1$. We will restrict our example to the case where the matrix has positive entries, it is possible to consider more general automorphisms under an invariant cone assumption, see \cite{AS}. Let $\Omega=\{0,1\}^{\N}$ and $\vartheta=\sigma$ be the left shift on $\Omega$. Let $A_0$, $A_1$ two hyperbolic automorphisms with positive entries. Let $A_0$ be chosen with a probability $q$ and $A_1$ with a probability $1-q$, i.e. $\P=P^{\N}$ with $P(0)=q$ and $P(1)=1-q$. Then the Lebesgue measure is stationary and the decay of correlations is exponentially fast for Lipschitz observables (indeed for strong H\"older observables (see \cite{AS})) $f$ and $g$ which satisfy $\int_{\mathbb{T}^2}f(x)dx=\int_{\mathbb{T}^2}g(x)dx=0$.
\begin{prop}\label{nonap}
This system is random-aperiodic.
\end{prop}
\begin{proof}
For any fixed $\omega\in\Omega$ and $n\in\N^*$, the matrix M corresponding to  $A^n_{\omega}=A_{\omega_n}\dots A_{\omega_2}A_{\omega_1}$ (i.e, $ A^n_{\omega}x=Mx \mod 1$) has positive integers entries. By the Perron-Frobenius Theorem, M has a simple largest eigenvalue $\lambda>0$ and since $|\det M|=1$, $\lambda>1$. Therefore, $A^n_{\omega}$ is also a hyperbolic toral automorphism. By \cite{PY}, we have for all $n\in\N^*$
\[\textrm{Card} \{x\in\mathbb{T}^2\,:\,A^n_{\omega}x=x\}<+\infty\]
and so
\[Leb(\{x\in\mathbb{T}^2\,:\,\exists n\in\mathbb{N}^*\,,\, A^n_{\omega}x=x\})=0.\]
Then, this system satisfies
\begin{equation*}
\P\otimes Leb(\{(\omega,x)\in\{0,1\}^{\N}\times\mathbb{T}^2\,:\,\exists n\in\N^*\,,\, A^n_{\omega}x=x\})=0.
\end{equation*}
\end{proof}
Proposition~\ref{nonap} and Corollary~\ref{cospdoc} give for $\P\otimes Leb$-almost every $(\omega,x)\in\{0,1\}^\N\times\mathbb{T}^2$
\[\lim_{r\rightarrow0}\frac{\log \tau_{r}^{\omega}(x)}{-\log r}=2\]
and Proposition~\ref{propegal} gives for $Leb$-almost every $x\in\mathbb{T}^2$
\[\lim_{r\rightarrow0}\frac{\log\left[\mathbf{T}_{r}(x)\right]}{-\log r}=2.\]
\section{Small random perturbations}
A particular case of random dynamical systems is the one where all the maps are chosen arbitraly close to a fixed initial map, namely the random dynamical system is a small random perturbation of a map $T$. More precisely we consider a measurable map $T:X\rightarrow X$  which plays the role of an initial dynamics, and we add a small amount of noise during its evolution, modelled by a random dynamical system. More precisely,
for any small $\ve>0$ (which is the noise level) we consider a probability space $(\Lambda_{\ve},\mathcal{B}(\Lambda_{\ve}),P_{\ve})$ where $\Lambda_{\ve}$ is a metric space and $\mathcal{B}(\Lambda_{\ve})$ the Borelian $\sigma$-algebra. Then let us consider a parametrized family of maps $\{T_{\lambda}\}_{\lambda\in \Lambda_{\ve}}$ which are $\ve$-close to $T$  in some $\mathcal{C}^0$ sense for small $\varepsilon$. The framework is the following $\Omega=\Lambda_{\eps}^\N$, $\P=P_{\eps}^\N$, $\vartheta=\sigma$ the left shift on $\Omega$ and $T_{\omega}=T_{\lambda_1}$ for all $\omega=(\lambda_1,\lambda_2,\dots)\in\Omega$. Let us consider the following map:
\begin{align*}
\Phi^{\ve}:~~\Lambda_{\ve}\times X&\longrightarrow X\\
(\lambda,x)&\longmapsto T_{\lambda}(x)
\end{align*}
such that for all $x\in X$, $\Phi_x^{\ve}:=\Phi^{\ve}(.,x) :\Lambda_{\ve}\rightarrow X$ is measurable and there is a $\lambda^{\ast}\in\Lambda_{\ve}$ which satisfies $\Phi^{\ve}(\lambda^{\ast},x)=T(x)$ for all $x\in X$.
We consider the probability measure $(\Phi_x^{\ve})_{\ast}P_{\ve}$ \textit{i.e}. the probability defined by:
$$
(\Phi_x^{\ve})_{\ast}P_{\ve}(A):=P_{\ve}\{\lambda\in\Lambda_{\ve}~:~T_{\lambda}x\in A\}
$$ 
Note that it is exactly the family of transition probabilities defining the associated Markov process of the random maps (see \ref{Markov}).
We set the  following classical  assumptions:
\begin{enumerate}[(RT1)]
\item For all $x\in X$, $(\Phi_x^{\ve})_{\ast}P_{\ve}$ is absolutely continuous with respect to the Lebesgue measure for any small $\varepsilon>0$.
\item For all $x\in X$,  $\Phi_x^{\ve}(\Lambda_{\ve})\subset B_{\ve}(Tx)$ for any small $\varepsilon>0$.
\end{enumerate}
\begin{remark} Our assumption  (RT2) implies a $\mathcal{C}^0$ closeness between $T$ and every $T_{\lambda}$.
\end{remark}
Again, we do not need non-instantaneous return time: 
\begin{prop}\label{propperiod}Under the hypothesis (RT1), the system is random-aperiodic.
\end{prop}
\begin{proof}
We will prove that for all $x\in X$
\[P_\eps^\N(\{\omega=(\lambda_1,\lambda_2,...)\in\Lambda_\eps^\N\,:\,\exists n\in\N^*\,,\, T^n_{\omega}x=x\})=0\]
where $T_{\omega}^n=T_{\lambda_n}\circ\dots\circ T_{\lambda_2}\circ T_{\lambda_1}$.

Let $x\in X$. Let $A\subset X$ with $Leb(A)=0$. Since the measure $(\Phi_x^{\ve})_{\ast}P_{\ve}$ is absolutely continuous by (RT1) we have $(\Phi_x^{\ve})_{\ast}P_{\ve}(A)=P_\eps(\{\lambda_1\in\Lambda_\eps\,:\,T_{\lambda_1}x\in A\})=0$. Moreover, we remark that
\[(\Phi_x^{\ve})^2_{\phantom{2}\ast}P_{\ve}^2(A):=P_\eps^2(\{(\lambda_1,\lambda_2)\in\Lambda_\eps^2\,:\,T_{\lambda_2}T_{\lambda_1}x\in A\})=0.\]
Indeed
\[(\Phi_x^{\ve})^2_{\phantom{2}\ast}P_{\ve}^2(A)=\int_{\Lambda_\eps}(\Phi_{T_{\lambda_1}x}^{\ve})_{\ast}P_{\ve}(A)dP_\eps(\lambda_1)\]
but since $T_{\lambda_1}x\in X$ for all $\lambda_1\in\Lambda_\eps$, (RT1) gives $(\Phi_{T_{\lambda_1}x}^{\ve})_{\ast}P_{\ve}(A)=0$ for all  $\lambda_1\in\Lambda_\eps$. Using this idea, with an easy induction argument, we can prove that for every $n\in\N^*$
\begin{equation}\label{phinpn}
(\Phi_x^{\ve})^n_{\phantom{n}\ast}P_{\ve}^n(A):=P_\eps^n(\{(\lambda_1,\dots,\lambda_n)\in\Lambda_\eps^n\,:\,T_{\lambda_n}\dots T_{\lambda_1}x\in A\})=0.
\end{equation}
Finally, we get
\begin{eqnarray*}
 P_\eps^\N(\{\omega\in\Lambda_\eps^\N:\exists n\in\N^*\,,\, T^n_{\omega}x=x\})&=& P_\eps^\N(\bigcup_{n\in\N^*}\{\omega\in\Lambda_\eps^\N\,:\, T^n_{\omega}x=x\})\\
 &\leq&\sum_{n\in\N^*} P_\eps^\N(\{\omega\in\Lambda_\eps^\N\,:\, T^n_{\omega}x=x\})\\
 &\leq&\sum_{n\in\N^*} P_\eps^n(\{(\lambda_1,\dots,\lambda_n)\in\Lambda_\eps^n\,:\, T_{\lambda_n}\dots T_{\lambda_1}x=x\})\\
 &\leq&\sum_{n\in\N^*} \underset{=0\textrm{  by \eqref{phinpn}}}{\underbrace{P_\eps^n(\{(\lambda_1,\dots,\lambda_n)\in\Lambda_\eps^n: T_{\lambda_n}\dots T_{\lambda_1}x\in\{x\}\})}}\\
 &\leq&0.
 \end{eqnarray*}
\end{proof}
\subsection{Small random perturbations of expanding maps of the circle}
Let $X=\mathcal{S}^1$ and $T$ be an expanding $\mathcal{C}^r$ $(2\leq r<\infty)$ transformation of $\mathcal{S}^1$. We put an additive noise to this system, namely  $T_{\lambda}x=Tx+\lambda$ where $\lambda$ is a random variable distributed  according to a density supported on $(-\ve,+\ve)$. Then it is well known (see for example \cite{BY}) that the random dynamical system admits an absolutely continuous stationary measure $\mu$ whose density is $\mathcal{C}^{r-1}$ and is exponentially mixing for $\mathcal{C}^{r-1}$ observables. Even if the decay of correlations is exponential for not Lipschitz observables, it is possible to go from $\mathcal{C}^{r-1}$ observables to Lipschitz observables with some simple approximation arguments. Then by Corollary~\ref{cospdoc} and Proposition~\ref{propperiod} for $\P\otimes\mu$-almost every $(\omega,x)\in(-\eps,\eps)^\N\times \mathcal{S}^1$ we have
\[\lim_{r\rightarrow0}\frac{\log \tau_{r}^{\omega}(x)}{-\log r}=1\]
and by Proposition~\ref{propegal} we get for $\mu$-almost every $x\in\mathcal{S}^1$
\[\lim_{r\rightarrow0}\frac{\log\left[\mathbf{T}_{r}(x)\right]}{-\log r}=1.\]
\subsection{Small random perturbations of  piecewise expanding maps of the interval}
Let $X=[0,1]$ and $T$ be a $\mathcal{C}^2$-piecewise expanding map without periodic turning point (see \cite{BY}) \footnote{Remark that similar results can be obtained for maps whose derivative is uniformly larger than 2, see \cite{KL}.}, then for the same additive perturbation than in the previous example, Baladi and Young \cite{BY} have proved the existence (and the stability) of an absolutely continuous stationary measure $\mu$ and the exponential decay of correlations for observables that are of bounded variations, and therefore for Lipschitz ones. We obtain by Corollary~\ref{cospdoc} and Proposition~\ref{propperiod} that for $\P\otimes\mu$-almost every $(\omega,x)\in(-\eps,\eps)^\N\times [0,1]$
\[\lim_{r\rightarrow0}\frac{\log \tau_{r}^{\omega}(x)}{-\log r}=1\]
and we obtain by Proposition~\ref{propegal} that for $\mu$-almost every $x\in[0,1]$
\[\lim_{r\rightarrow0}\frac{\log\left[\mathbf{T}_{r}(x)\right]}{-\log r}=1.\]


\section{Proofs}\label{proof}
We will use the notion of   weakly diametrically regular measures:
\begin{defi}
A measure $\mu$ is weakly diametrically regular (wdr) on the set $Z\subset X$ if for any $\eta>1$, for $\mu$-almost every $x\in Z$ and every $\eps>0$, there exists $\delta>0$ such that if $r<\delta$ then $\mu\left(B\left(x,\eta r\right)\right)\leq\mu\left(B\left(x,r\right)\right)r^{-\eps}$.
\end{defi}
%
%
\begin{proof}[Proof of Theorem~\ref{thborne}]
It is well known (see for instance \cite{BS}) that   any probability measure is weakly diametrically regular on $\R^d$  for all $d\in\N^*$ and then so is the measure $\mu$ on $X$. Let us emphasize that in the previous definition   the function $\delta(\cdot,\eps,\eta)$ can be made measurable for every fixed $\eps$ and $\eta$. Let us fix $\eps>0$ and $\eta=4$. We choose $\delta>0$   small enough to get: $$\mu(X_\delta)>\mu(X)-\eps=1-\eps$$
 where $X_\delta:=\left\{x\in X:\,\delta(x,\eps,\eta)>\delta\right\}.$
We now set the following lemma in order to use the Borel-Cantelli's one:
\begin{lem}
For every $p\in\N$, $\displaystyle{\sum_{n\in\N}\mu(A_\eps(e^{-n}))<\infty}$ where for $r>0$
\[A_\eps(r):=\left\{y\in X_\delta: \mathbf{T}_{6r,p}(y)\mu\left(B\left(y,2r\right)\right)\geq r^{-2\eps}\right\}.\]
\end{lem}

\begin{proof}
\begin{defi} Given $r>0$, a countable set $E\subset F$ is a maximal $r$-separated set for $F$  if
\begin{enumerate}
\item $B(x,\frac{r}{2})\cap B(y,\frac{r}{2})=\emptyset$ for any two distinct $x,y\in E$.
\item $\mu(F\,\backslash \underset{x\in E}{\bigcup}B(x,r))=0$.
\end{enumerate}
\end{defi}
Let $p\in\N$, $r>0$ and $C\subset X_\delta$ be a maximal $2r$-separated set for $X_\delta$.
\begin{eqnarray}
\mu(A_\eps(r))&=&\mu\left(\left\{y\in X_\delta: \mathbf{T}_{6r,p}(y)\mu\left(B\left(y,2r\right)\right)\geq r^{-2\eps}\right\}\right)\nonumber\\
&\leq&\sum_{x\in C}\mu\left(\left\{y\in B\left(x,2r\right):\,\mathbf{T}_{6r,p}(y)\mu\left(B\left(y,2r\right)\right)\geq r^{-2\eps}\right\}\right).\label{ensmaxsep}
\end{eqnarray}
For $y\in B(x,4r)$, we define: \[\tau_{4r,p}^{\omega}(y,x):=\inf\left\{k>p:\,d(T_{\vartheta^k\omega}\circ\dots\circ T_{\omega}y,x)<4r\right\}\] and:  \[\mathbf{T}_{4r,p}(y,x):=\int_{\Omega}\tau_{4r,p}^{\omega}(y,x)d\P(\omega).\] 
If $d(x,y)<2r$, for all $\omega\in\Omega$ we have
$\tau^{\omega}_{4r,p}(y,x)\geq\tau^{\omega}_{6r,p}(y)$
and then:
\begin{equation}\label{inegarx}
\mathbf{T}_{4r,p}(y,x)\mu\left(B\left(x,4r\right)\right)\geq\mathbf{T}_{6r,p}(y)\mu\left(B\left(y,2r\right)\right).
\end{equation}
It follows that for any $x\in C$:
\begin{equation}\label{drx}
\mu\left(\left\{y\in B\left(x,2r\right):\,\mathbf{T}_{6r,p}(y)\mu\left(B\left(y,2r\right)\right)\geq r^{-2\eps}\right\}\right)\leq \mu(D_{r,x})
\end{equation}
where:
\[D_{r,x}:=\left\{y\in B\left(x,4r\right):\,\mathbf{T}_{4r,p}(y,x)\mu\left(B\left(x,4r\right)\right)\geq r^{-2\eps}\right\}.\]
We then use the Markov's inequality and we get:
\begin{eqnarray}\label{cheb}
\mu(D_{r,x})&\leq& r^{2\eps}\mu\left(B\left(x,4r\right)\right)\int_{B(x,4r)}\mathbf{T}_{4r,p}(y,x)\,d\mu(y)\\
&=& r^{2\eps}\mu\left(B\left(x,4r\right)\right)\int_{B(x,4r)\times\Omega}\tau^{\omega}_{4r,p}(y,x)\,d\mu(y)d\P(\omega).
\end{eqnarray}
Since $\tau^{\omega}_{4r,p}(y,x)$ is bounded by the $p^{th}$ return time of $(y,\omega)$ in the set $B(x,4r)\times\Omega$, Kac's lemma provides the following inequality:
\begin{equation}\label{kacbosh}
\int_{B(x,4r)\times\Omega}\tau^{\omega}_{4r,p}(y,x)\,d\mu(y)d\P(\omega)\leq p.
\end{equation}
Combining (\ref{cheb}) with (\ref{kacbosh}), we have:
\begin{equation}\label{ineg}
\mu(D_{r,x})\leq pr^{2\eps}\mu\left(B\left(x,4r\right)\right)
\end{equation}
and thus:
\begin{eqnarray*}
\mu(A_\eps(r))&\leq&\sum_{x\in C}\mu(D_{r,x}) \qquad \textrm{by \eqref{ensmaxsep} and \eqref{drx}}\\
&\leq&p\,r^{2\eps}\sum_{x\in C}\mu\left(B(x,4r)\right)  \qquad \textrm{by (\ref{ineg})}\\
&\leq&p\,r^{\eps}\sum_{x\in C}\mu\left(B(x,r)\right)  \qquad \textrm{since $\mu$ is wdr}\\
&\leq&p\,r^\eps \qquad  \textrm{by definition of $C$}.
\end{eqnarray*}
Finally:
\[\sum_{n,e^{-n}<\delta}\mu(A_\eps(e^{-n}))=\sum_{n>-\log\delta}\mu(A_\eps(e^{-n}))\leq p\sum_ne^{-\eps n}<\infty.\]
\end{proof}
We can thus apply the Borel-Cantelli lemma to get for large $n$ and for $\mu$-almost every $x\in X_\delta$
\[\mathbf{T}_{6e^{-n},p}(x)\mu\left(B(x,2e^{-n})\right)\leq e^{2\eps n}\]
\textit{i.e.}
\begin{equation}\label{inegfin}
\frac{\log\mathbf{T}_{6e^{-n},p}(x)}{n}\leq 2\eps+\frac{\log\mu(B(x,2e^{-n}))}{-n}.
\end{equation}
One can easily prove that for all $a>0$ we have:
\begin{eqnarray*}
\underline{d}_\mu(x)=\underset{n\rightarrow\infty}{\liminf}\frac{\log\mu\left(B\left(x,ae^{-n}\right)\right)}{-n}&\textrm{ and }&\overline{d}_\mu(x)=\underset{n\rightarrow\infty}{\limsup}\frac{\log\mu\left(B\left(x,ae^{-n}\right)\right)}{-n}\\
\underline{\mathbf{R}}(x)=\lim_{p\rightarrow\infty}\liminf_{n\rightarrow\infty}\frac{\log\mathbf{T}_{ae^{-n},p}(x)}{n}&\textrm{ and }&\overline{\mathbf{R}}(x)=\lim_{p\rightarrow\infty}\limsup_{n\rightarrow\infty}\frac{\log\mathbf{T}_{ae^{-n},p}(x)}{n}
\end{eqnarray*}
Since $\eps$ can be chosen arbitrarily small we get the result considering the inferior (resp. the superior) limit when $n\rightarrow+\infty$ and the limit when $p\rightarrow+\infty$  in (\ref{inegfin}).
\end{proof}


\begin{proof}[Proof of Theorem~\ref{thegal}]
Let $x\in X$ be such that $\underline{R}^{\omega}(x)\geq \underline{d}_\mu(x)>0$ for $\P$-almost every $\omega\in\Omega$, the existence of such a $x$ is ensured by Theorem~\ref{spdoc}. Let  $0<\eps<1$, by the proof of Theorem 5 in \cite{RS}, we know that for $\P$-almost every $\omega\in\Omega$ it exists $N\in\N$ such that for every $p>N$, $\underset{r\rightarrow0}{\liminf}\,\mu(B(x,r))^{1-\eps}\tau_{r,p}^{\omega}(x)=+\infty$, that is to say that:
\begin{equation}
\P(\left\{\omega\in\Omega:\,\exists N\in\N,\,\forall p>N,\,\underset{r\rightarrow0}{\liminf}\,\mu(B(x,r))^{1-\eps}\tau_{r,p}^{\omega}(x)=+\infty\right\})=\P(\Omega)=1.
\end{equation}
Let us denote by 
$$\tilde\Omega(N):=\left\{\omega\in\Omega\,:\,\forall p>N,\,\underset{r\rightarrow0}{\liminf}\,\mu(B(x,r))^{1-\eps}\tau_{r,p}^{\omega}(x)=+\infty\right\}.$$
It exists $N_1\in\N$ such that $\P(\tilde\Omega(N_1))>0 $, otherwise $\P(\cup_{N\in\N}\tilde\Omega(N))\leq\sum_{N\in\N}\P(\tilde\Omega(N))=0$ and it would contradict the fact that:
$$
1=\P(\bigcup_{N\in\N}\tilde\Omega(N))=\P(\left\{\omega\in\Omega:\,\exists N\in\N,\,\forall p>N,\,\underset{r\rightarrow0}{\liminf}\,\mu(B(x,r))^{1-\eps}\tau_{r,p}^{\omega}(x)=+\infty\right\}).
$$
Let $p>N_1$,  we therefore have:
\begin{equation}\label{intliminfini}
\int_{\tilde\Omega}\underset{r\rightarrow0}{\liminf}\,\mu(B(x,r))^{1-\eps}\tau_{r,p}^{\omega}(x)d\P(\omega)=+\infty.\end{equation}
Moreover, Fatou's lemma gives
\begin{eqnarray*}
+\infty=\int_{\tilde\Omega}\underset{r\rightarrow0}{\liminf}\,\mu(B(x,r))^{1-\eps}\tau_{r,p}^{\omega}(x)d\P(\omega)&\leq&\int_{\Omega}\underset{r\rightarrow0}{\liminf}\,\mu(B(x,r))^{1-\eps}\tau_{r,p}^{\omega}(x)d\P(\omega)\\&\leq&\underset{r\rightarrow0}{\liminf}\int_{\Omega} \mu(B(x,r))^{1-\eps}\tau_{r,p}^{\omega}(x)d\P(\omega)\nonumber\\
&=&\underset{r\rightarrow0}{\liminf}\,\mu(B(x,r))^{1-\eps}\mathbf{T}_{r,p}(x).
\end{eqnarray*}
We thus get that:
\[\underset{r\rightarrow0}{\liminf}\,\mu(B(x,r))^{1-\eps}\mathbf{T}_{r,p}(x)=\underset{r\rightarrow0}{\limsup}\,\mu(B(x,r))^{1-\eps}\mathbf{T}_{r,p}(x)=+\infty.\]
Finally we have that for all $p>N_1$ and $M>0$, it exists $R>0$ such that for all $r<R$, $\mu(B(x,r))^{1-\eps}\mathbf{T}_{r,p}(x)\geq M$ and then:
$$
(1-\eps)\frac{\log\mu(B(x,r))}{\log r}\leq \frac{M}{\log r}-\frac{\log\mathbf{T}_{r,p}(x)}{\log r}
$$
which drives to:
\[\underline{\mathbf{R}}(x)\geq (1-\eps)\underline{d}_\mu(x)\qquad\textrm{and}\qquad\overline{\mathbf{R}}(x)\geq (1-\eps)\overline{d}_\mu(x).\]
Since these inequalities hold for every $0<\eps<1$, we get:
\begin{equation}\label{eqrbigo}
\underline{\mathbf{R}}(x)\geq \underline{d}_\mu(x)\qquad\textrm{and}\qquad\overline{\mathbf{R}}(x)\geq \overline{d}_\mu(x).\end{equation}
By Theorem~\ref{spdoc}, the equation \eqref{eqrbigo} is satisfied for $\mu$-almost every $x$ such that $\underline{d}_\mu(x)>0$ and then the theorem is proved using Theorem~\ref{thborne}.
\end{proof}
As in the previous proof, the principal idea of the proof of Proposition~\ref{propegal} is to use Fatou's lemma:
\begin{proof}[Proof of Proposition~\ref{propegal}]
Let $x\in X$ be such that $$\underset{r\rightarrow0}{\liminf} \frac{\log\tau_r^\omega(x)}{-\log r}\geq \underline{d}_\mu(x)>0$$ for $\P$-almost every $\omega\in\Omega$; the existence of such a $x$ is ensured by Corollary~\ref{cospdoc}. Let  $0<\eps<1$, for $\P$-almost every $\omega\in\Omega$ we have $\underset{r\rightarrow0}{\liminf}\,\mu(B(x,r))^{1-\eps}\tau_{r}^{\omega}(x)=+\infty$ and then
\begin{equation}\label{intliminfini2}
\int_{\Omega}\underset{r\rightarrow0}{\liminf}\,\mu(B(x,r))^{1-\eps}\tau_{r}^{\omega}(x)d\P(\omega)=+\infty.\end{equation}
By Fatou's lemma we have
\begin{eqnarray*}
+\infty=\int_{\Omega}\underset{r\rightarrow0}{\liminf}\,\mu(B(x,r))^{1-\eps}\tau_{r}^{\omega}(x)d\P(\omega)&\leq&\underset{r\rightarrow0}{\liminf}\int_{\Omega} \mu(B(x,r))^{1-\eps}\tau_{r}^{\omega}(x)d\P(\omega)\nonumber\\
&=&\underset{r\rightarrow0}{\liminf}\,\mu(B(x,r))^{1-\eps}\mathbf{T}_{r}(x).
\end{eqnarray*}
So we obtain
\[\underset{r\rightarrow0}{\liminf}\,\mu(B(x,r))^{1-\eps}\mathbf{T}_{r}(x)=\underset{r\rightarrow0}{\limsup}\,\mu(B(x,r))^{1-\eps}\mathbf{T}_{r}(x)=+\infty.\]
Then, it exists $R>0$ such that for all $r<R$ we have $\mu(B(x,r))^{1-\eps}\mathbf{T}_r(x)\geq M$ and so:
$$
(1-\eps)\frac{\log\mu(B(x,r))}{\log r}\leq \frac{M}{\log r}-\frac{\log\mathbf{T}_r(x)}{\log r}
$$
which gives us:
\[\liminf_{r\rightarrow0}\frac{\log\mathbf{T}_r(x)}{\log r}\geq (1-\eps)\underline{d}_\mu(x)\qquad\textrm{and}\qquad\limsup_{r\rightarrow0}\frac{\log\mathbf{T}_r(x)}{\log r}\geq (1-\eps)\overline{d}_\mu(x).\]
These inequalities are satisfied for $\eps$ arbitrarily small, so we get:
\begin{equation}\label{eqrbigo2}
\liminf_{r\rightarrow0}\frac{\log\mathbf{T}_r(x)}{\log r}\geq \underline{d}_\mu(x)\qquad\textrm{and}\qquad\limsup_{r\rightarrow0}\frac{\log\mathbf{T}_r(x)}{\log r}\geq\overline{d}_\mu(x).\end{equation}
Since \eqref{eqrbigo2} is satisfied for $\mu$-almost every $x$ such that $\underline{d}_\mu(x)>0$ by Corollary~\ref{cospdoc}, we get the result by Corollary~\ref{coborne}.
\end{proof}

{\small\textit{Acknowledgement.} The authors would like to thank B. Saussol for his useful help and comments.}

\bibliographystyle{siam} 
\bibliography{biblio}

\end{document}